\documentclass[letterpaper, 10 pt, conference]{ieeeconf}  

\IEEEoverridecommandlockouts    

\usepackage{graphics} 
\usepackage{amsmath} 
\usepackage{amssymb}  
\usepackage{subfigure}
\usepackage{stfloats}
\usepackage[T1]{fontenc}
\usepackage{graphicx}
\usepackage{algorithm,algorithmic}

\newcommand{\X}{\mathbb{X}}
\newcommand{\R}{\mathbb{R}}

\newcommand{\C}{\mathbb{C}}

\newcommand{\K}{\mathcal{K}}
\newcommand{\PP}{\mathbb{P}}
\newcommand{\E}{\mathbb{E}}
\newcommand{\gen}{\mathcal{L}}
\newcommand{\A}{\mathcal{A}}
\newcommand{\Abs}[1]{\left|#1\right|}
\newcommand{\Par}[1]{\left(#1\right)}
\newcommand{\norm}[1]{\left\|#1\right\|}
\newcommand{\FP}{\phi_\lambda}
\newcommand{\FPi}{\phi_{\lambda_i}}
\usepackage[dvipsnames]{xcolor}
\usepackage[space]{cite}

\newtheorem{defn}{Definition}
\newtheorem{prop}{Proposition}
\newtheorem{thm}{Theorem}

\newtheorem{ex}{Example}
\newtheorem{remark}{Remark}%

\title{\LARGE \bf Stability analysis of stochastic differential equations \\ with the Koopman operator
}

\author{François-Grégoire Bierwart$^{1}$, Gaëtan Louvet$^{2}$ and Alexandre Mauroy$^{3}$
\thanks{$^{1,3}$ Department of Mathematics and Namur Research Institute for Complex Systems (naXys), University of Namur, 5000, Belgium
        {\tt\small francois-gregoire.bierwart.@unamur.be, alexandre.mauroy@unamur.be}}
\thanks{$^{2}$ Université libre de Bruxelles (ULB),
        {\tt\small gaetan.louvet@ulb.be}}%
}

\begin{document}

\maketitle
\thispagestyle{empty}
\pagestyle{empty}

\begin{abstract}
This paper investigates the stability of stochastic systems through the lens of the Koopman operator framework. In particular, we establish the connection between the convergence properties of the Koopman semigroup and specific stability properties of the underlying stochastic system. Moreover, we demonstrate that specific Koopman eigenfunctions capture stochastic attractivity properties and we illustrate these results with two numerical examples. This work mirrors existing Koopman operator-based stability results for deterministic systems and extends them to the stochastic framework.
\end{abstract}

\section{Introduction}
Stochastic systems are crucial for modeling many real-life applications. Also, the analysis of their long-term behavior is of paramount importance and, to that end, several notions of stability can be considered, including the stochastic counterpart of Lyapunov stability (see e.g. \cite[Chapter 4]{mao2007stochastic}), which can be captured via the design of Lyapunov functions \cite{mao2007stochastic,khasminskii2011stochastic,mao1999stochastic}.  In this context, several Lyapunov function design methods have been developed over the last decades, such as sum-of-squares methods (SOS) \cite{hafstein2018lyapunov}, solutions to stochastic Zubov equation \cite{CamilliGrune2003}, and meshless collocation-based methods \cite{bjornsson2018computation,bjornsson2018verification}. Moreover, recent advances in machine learning also contributed to the development of novel Lyapunov design techniques, see e.g. \cite{lechner2022stability,su2025verifying} and references therein.

An alternative characterization of stability properties is offered by the Koopman operator framework. The Koopman operator is an infinite-dimensional operator that provides a linear description of the original
nonlinear system in terms of the evolution of the expected value of observables. Moreover, a Lyapunov function can be seen as a specific observable that decreases under the action of the Koopman operator. In the case of deterministic systems, the work \cite{mauroy2014global} established a relationship between stability properties and spectral properties of the associated Koopman operator, showing in particular that a generic Lyapunov function can be constructed directly from Koopman eigenfunctions. However, to the best of our knowledge, there is no such work that assesses the stability of stochastic systems through the spectral properties of the Koopman operator. We will fill this gap in the present paper.

This paper builds upon the work \cite{mauroy2014global} by extending several of its results to stochastic systems. More precisely, we first establish a relationship between the stability properties of the system and the asymptotic behavior of the Koopman operator applied to specific observables. Due to the limited applicability of this result, we investigate the spectral properties of the Koopman operator and we show that the system stability is captured by specific Koopman eigenfunctions.       

The paper is organized as follows. In Section \ref{sec:preli}, we present the fundamentals of stochastic differential equations and related notions of stability. 
The connection between the Koopman operator convergence properties and the stability properties is investigated in Section \ref{sec:preli_res}. In Section \ref{sec:res_eigen}, we develop our main results relying on the spectral properties of the Koopman operator. These results are further illustrated with two examples in Section \ref{sec:numerical}. Finally, concluding remarks and perspectives are given in Section \ref{sec:ccl}.

\section{Preliminaries}\label{sec:preli}

\noindent In this section, we present the basics of stochastic differential equations and associated stability properties.   

\subsection{A brief introduction to stochastic differential equations}

We introduce the notion of stochastic differential equations with the following definitions. 
\begin{defn}
Let $(\Omega,\mathcal{P},\mathbb{P})$ be a probability space. A filtration is a family $\{\mathcal{P}_t\}_t$ of increasing sub$-\sigma-$algebra of $\mathcal{P}$ (i.e. $\mathcal{P}_t\subset\mathcal{P}_s\subset\mathcal{P}$ for all $0\leq t < s < \infty$). The filtration is said to be right continuous if $\mathcal{P}_t = \cap_{s>t}\mathcal{P}_s$  for all $t\geq 0$ . When the probability space is complete, the filtration is said to satisfy the usual conditions if it is right continuous and $\mathcal{P}_0$ contains all $\mathbb{P}-$null sets. 
\end{defn}
\begin{defn}
Let $(\Omega,\mathcal{P},\mathbb{P})$ be a probability space and $\{\mathcal{P}_t\}_t$ a filtration. A stochastic process $\{Y_t\}_t$ is adapted to the filtration $\mathcal{P}_t$ if, for each $t\geq 0$, $Y_t$ is $\mathcal{P}_t-$measurable.  
\end{defn}
\begin{defn}
Let $(\Omega,\mathcal{P},\mathbb{P})$ be a probability space and $\{\mathcal{P}_t\}_t$ a filtration. For any $T_1,T_2\geq0$, let $\mathcal{M}^2([T_1,T_2],\R^n)$ be the set of all $\{\mathcal{P}_t\}-$adapted $\R^n$ valued processes $\{Y_t\}_t$ such that 
\begin{equation}\label{eq:M2}
\mathbb{E}\left[\int_{T_1}^{T_2} \|Y_t\|^2\,\mathrm{d}t\right]<\infty
\end{equation}
where $\mathbb{E}[\cdot]$ denotes the expected value.  
\end{defn}
\noindent Let $W_t = (W_1(t),W_2(t),\ldots,W_m(t))$ be a $m-$dimensional Wiener process on the complete probability space $(\Omega,\mathcal{P},\mathbb{P})$. We define the filtration $\mathcal{P}_t = \sigma\left(\left\{W_s^{-1}(B) \mid 0\leq s \leq t, B \in\mathcal{B}_{\R}\right\}\right)$ augmented with the $\mathbb{P}-$null sets, where for any $A\subset\Omega$, $\sigma(A)$ is the smallest $\sigma-$algebra that contains $A$.  We make this choice so that $\{W_t\}_t$ is adapted to $\{\mathcal{P}_t\}_t$ and $\{\mathcal{P}_t\}_t$ satisfies usual conditions.

\noindent An autonomous stochastic differential equation (SDE) is an equation of the form  
\begin{equation}\label{SDE}
\left\{
\begin{array}{l}
dX_t = F(X_t)dt + G(X_t)dW_t,\\[0.2cm]
X_0 = x\in \R^n, 
\end{array}
\right.\vspace{0.2cm}
\end{equation}
where $F :\R^n \rightarrow \R^n, G : \R^n\rightarrow \R^{n\times m}$ and 
$
W_t = (W_1(t),W_2(t),\ldots,W_m(t))
$
is an $m-$dimensional Wiener process. This form is a short-hand for the integral equation 
\begin{equation}\label{Xt}
X_t^{x} = x + \int_0^tF(X_s^x)ds + \int_0^tG(X_s^x)dWs, \quad t \geq 0, 
\end{equation}
where the second integral is an Itô integral and $\{X_t^x\}_{t\geq0}$ is a stochastic process where, for any $t\geq 0$, $X_t^x : \Omega \rightarrow \R^n$. We refer the reader to \cite{mao2007stochastic,oksendal2013stochastic} for more details about SDEs and Itô integrals. The following theorem \cite[Theorem 3.5] {mao2007stochastic} establishes the existence and uniqueness of the solution to \eqref{Xt}. 

\begin{thm}\label{thm:solSDE}
For any real numbers $a,b$, let $a \vee b := \max(a,b)$. Assume that the following holds: 
\begin{enumerate}
    \item (Local Lipschitz condition) For every integer $k\geq 1$, there exists a positive constant $\overline{C}_k$ such that, for all $x,y\in\R^n$ with $\|x\|\vee\|y\|<k$, 
    $$
    \|F(x)-F(y)\|^2 \;\bigvee\; \|G(x)-G(y)\|^2 \leq \overline{C}_n\|x-y\|^2;
    $$
    \item (Monotonicity condition) There exists a positive constant $C$ such that, for all $x\in\R^n$, 
    $$
    x^\top F(x) + \dfrac{1}{2}\mathrm{Tr}(G(x)^\top G(x)) \leq C(1+\|x\|^2)
    $$
    where $\mathrm{Tr}(\cdot)$ denotes the trace. 
\end{enumerate}
Then, there exists a unique solution to \eqref{Xt} which belongs to $\mathcal{M}^2([0,\infty[,\R^n)$. This solution is referred to as a strong solution.
\end{thm}
\begin{remark}
Standard conditions used in the literature to establish the existence and uniqueness of solutions to \eqref{Xt} are the global Lipschitz condition along with the linear growth condition. However, these conditions are more restrictive than those of Theorem \ref{thm:solSDE}. In particular, no polynomial of degree greater than two 
satisfies both conditions, while the assumptions of the theorem may hold for some polynomial functions, as illustrated in Example \ref{ex3}. 
\end{remark}

\subsection{Stability of stochastic differential equations}

This section presents several notions of stability for SDEs, which are adapted from classic definitions related to deterministic systems. We recall these definitions below. 

\begin{defn}[Stability for ODEs]\label{def:stablility:ODEs}
Consider the system \eqref{SDE} with $G = 0$ and let $A\subset \R^n$ be an invariant set. For any $z\in \R^n$, we define the distance
$$
d(z,A) = \displaystyle\inf_{y\in A}\|z-y\|,
$$
between $z$ and $A$. 
Then, 

\begin{enumerate}
    \item $A$ is stable if, for all $\varepsilon>0$, there exists $\delta > 0$ such that 
    $$
    d(x,A)<\delta \Longrightarrow \sup_{t\geq 0}\,d(X_t^x,A)<\varepsilon.
    $$
    \item $A$ is pointwise attractive in $D\subset \R^n$ if, for all $x\in D$,
    $$
    \lim_{t\rightarrow\infty}d(X_t^x,A) = 0.\\[0.3cm]
    $$
\end{enumerate}
\end{defn}

\noindent For SDEs, each initial condition $x\in\X$ is associated with an infinite collection of trajectories $\{X_t^x(p)\}_{p\,\in\,\Omega}$, corresponding to different realizations. As a result, different dynamical behaviors may occur, leading to several notions of attractiveness \cite{khasminskii2011stochastic,mao2007stochastic}.

\begin{defn}[Attractivity for SDEs]\label{def:attractor:sde}
Let $A\subset \R^n$. Then,
\begin{enumerate}
\setlength{\itemsep}{0.2cm}
    \item[1.] $A$ is stable in probability if, for all $r > 0$ and $\varepsilon \in (0,1)$, there exists $\delta > 0$ such that 
    $$
    d(x,A) < \delta \Longrightarrow \PP\left(\sup_{t\geq0}\,d(X_t^x,A)<r\right) \geq 1-\varepsilon.
    $$
    \item[2.] $A$ is attractive in probability if, for all $\varepsilon\in(0,1)$, there exists $\delta>0$ such that 
    $$
    d(x,A)<\delta \Longrightarrow \PP\left(\displaystyle\lim_{t\rightarrow\infty}\,d(X_t^x,A) = 0\right)\geq 1-\epsilon
    $$
    \item[3.] $A$ is  asymptotically stable in probability if it is stable and attractive in probability. 
    \item[4.] $A$ is attractive in mean in $D\subset\R^n$ if, for all $x \in D$,
    $$
    \lim_{t\rightarrow\infty}\,\E\left[d(X_t^x,A)\right] = 0.
    $$
    \item[5.] $A$ is almost surely attractive in $D\subset\R^n$ if, for all $x \in D$, 
    $$
    \PP\left(\lim_{t\rightarrow\infty}\,d(X_t^x,A) = 0\right) = 1,
    $$
    that is $\displaystyle\lim_{t\rightarrow\infty}d(X_t^x,A) = 0$ almost surely (a.s.) for all $x\in D$.
\end{enumerate}
\end{defn}

\section{Stability results with the Koopman operator}\label{sec:preli_res}

In this section, we introduce the Koopman semigroup associated with stochastic systems and investigate its connections to the stability of SDEs.

\subsection{Koopman operator framework}

From this point on, we suppose that the assumptions of Theorem \ref{thm:solSDE} are satisfied. Moreover, we assume that $\X\subset\R^n$ is almost surely  forward invariant, that is, $X_t^x(p) \in \X$ for all $x\in\X$, for all $t\geq 0$, and for almost every $p\in\Omega$.

For stochastic systems, the Koopman operator is defined as follows \cite[Chapter 6]{mauroy2020koopman} : 
\begin{defn}\label{KO}
Let $\mathcal{F}$ be a Banach space of observables \mbox{$f : \X\rightarrow\C$}. The Koopman semigroup $(K_t)_{t\geq0}$ with $K_t : \mathcal{F} \rightarrow \mathcal{F}$ is defined by 
$$(K_tf)(x) = \mathbb{E}\left[f(X_t^x)\right] $$
for any $f\in\mathcal{F}$ and $x\in\X$. 
\end{defn}
In \cite{mauroy2014global}, the authors establish a connection between the global attractivity properties of deterministic systems and the convergence properties of the associated Koopman semigroup. For stochastic systems, the classic definition of Koopman semigroup (Definition \ref{KO}) is not sufficient to capture the various definitions of stability. This motivates the additional definition of the pathwise Koopman semigroup.
\begin{defn}\label{KO_path}
Let $\mathcal{F}$ be a Banach space of observables \mbox{$f : \X\rightarrow\C$.} For any $p\in\Omega$, the pathwise Koopman semigroup $(K_{t,p})_{t \geq0}$ with $K_{t,p}\ : \mathcal{F} \rightarrow \mathcal{F}$ is defined by 
$$(K_{t,p}f)(x) = f(X_t^x(p))$$
for any $f\in\mathcal{F}$ and $x\in\X$.    
\end{defn}
\begin{remark}
For any $f\in\mathcal{F}$ and $x\in\X$, we have that $
(K_t f)(x) = \mathbb{E}\left[(K_{t,\cdot}f)(x)\right]$.  
\end{remark}
Moreover, we introduce the following operators, which will be used for the spectral analysis in Section \ref{sec:res_eigen}.  
\begin{defn}
For any $x\in\X$, we define:  
\begin{itemize}
\setlength{\itemsep}{0.2cm}
    \item The Koopman generator 
    $$(\mathcal{A} f)(x) = \lim_{t\rightarrow 0^+} \dfrac{(K_tf)(x) - f(x)}{t}$$
    for every $f\in\mathcal{F}$ whenever the limit exists ; 
    \item The SDE generator 
    $$
    (\gen f)(x) = (\nabla f^\top F)(x) + \dfrac{1}{2}\mathrm{Tr}\left(G^\top\nabla^2f G\right)(x)
    $$
    for any $f\in C^2(\X)$. \vspace{0.2cm}
\end{itemize}
\end{defn}
Let $C^2_b(\X)$ be the space of twice continuously differentiable bounded functions on $\X$, with bounded first and second derivatives. It can be shown  that $\A f=\gen f$ for all $f\in C^2_b(\X)$  (see e.g. \cite[Chapter 6]{mauroy2020koopman}). 

\subsection{First stability results}

We are now in a position to extend the results presented in \cite{mauroy2014global} to the case of stochastic systems. For $A\subset\X$, let 
$$
\mathcal{F}_{A_c} = \left\{f\in\mathcal{F} \mid f(x) = 0\;\forall x \in A\right\}
$$
be the set of functions vanishing on $A$.

The following result shows the connection between the stability properties of the SDE and those of the associated Koopman semigroup. 
\begin{prop}\label{prop:1} Let the system \eqref{SDE} and $A\subset \X$ be compact. Then, 
\begin{enumerate}
\setlength{\itemsep}{0.2cm}
\item $A$ is attractive in mean in $D\subset\X$ if and only if 
$$\lim_{t\rightarrow\infty}(K_tf)(x) = 0 \quad\forall f\in\mathcal{F}_{A_c}, \forall x \in D,$$ 
where $\mathcal{F} = \mathrm{Lip}(\X)$ is the set of globally Lipschitz continuous functions on $\X$. 
\item $A$ is almost surely attractive in $D \subset \X$ if and only if, for all $f\in\mathcal{F}_{A_c}$ and $x \in D$, 
$$
\lim_{t\rightarrow\infty}(K_{t,p}f)(x) = 0 
$$
for almost every $p\in\Omega$ with $\mathcal{F} = C(\X)$. 
\item $A$ is stable in probability if and only if, for all $f\in\mathcal{F}_{A_c}$, $\varepsilon\in(0,1)$, and $r > 0$, there exists $\delta > 0$ such that 
$$
d(x,A)<\delta \Longrightarrow \PP\left(\sup_{t\geq0}\left|(K_{t,\cdot}f)(x)\right|<r\right)\geq 1-\varepsilon
$$
with $\mathcal{F} = C(\X)$.
\item $A$ is attractive in probability if and only if, for all $f\in\mathcal{F}_{A_c}$ and $\varepsilon\in(0,1)$, there exists $\delta > 0$ such that 
$$
d(x,A)<\delta \Longrightarrow \PP\left(\lim_{t\rightarrow\infty}\left|(K_{t,\cdot}f)(x)\right|=0\right)\geq 1-\varepsilon
$$
with $\mathcal{F} = C(\X)$.
\end{enumerate}
\end{prop}

\begin{proof}~
\begin{enumerate}
\setlength{\itemsep}{0.2cm}
\item  \textit{Sufficiency}. The result holds with $f(\cdot) = d(\cdot,A)$, which is globally Lipschitz continuous on $\R^n$ and has zero value on $A$.\\
\textit{Necessity}. Let $f\in\mathcal{F}_{A_c}$. Since $f$ is Lipchitz continuous on $\X$ and has zero value on $A$, there exists $L>0$ such that, for any $x\in D$ and any $a\in A$, 
\begin{equation}\label{eqq}
\Abs{f\Par{X_t^x}} = \Abs{f\Par{X_t^x}-f(a)} \leq L\norm{X_t^x-a}.
\end{equation}
Since \eqref{eqq} holds for any $a\in A$, it also holds for the infimum so that  
$$
\Abs{f\Par{X_t^x}} \leq L d(X_t^x,A).
$$
By monotonicity of the expected value together with the Jensen's inequality, we obtain
$$
\Abs{\E\left[f\Par{X_t^x}\right]} \leq \E\left[\Abs{f\Par{X_t^x}}\right]\leq L\E\left[d(X_t^x,A)\right].
$$
Since $A$ is mean attractive, we finally have
$$\displaystyle\lim_{t\rightarrow\infty}\E\left[f\Par{X_t^x}\right] = 0.$$ 

\item \textit{Sufficiency}. The result holds with $f(\cdot) = d(\cdot,A)$, which is continuous on $\R^n$ and has zero value on $A$. \\
\textit{Necessity}. We first recall that, for any family of measurable functions $(g_t)_{t\geq 0}$ on the complete measure space $\Omega$, $g_t$ converges almost surely to $0$ if and only if there exists
$N\subset \Omega$ with $\PP\left(N\right) = 0$ such that, for all $p\in N^c$ and $\varepsilon > 0$, there exists $T_{p,\varepsilon} > 0$ with $|g_t(p)| < \varepsilon$ for all $t\geq T_{p,\varepsilon}$. Let us fix $f\in\mathcal{F}_{A_c}$, $x\in D$, and $\varepsilon>0$. Since $f$ is continuous and has zero value on $A$, there exists $\delta_\varepsilon>0$ such that 
$$
d(x,A)<\delta_\varepsilon\Longrightarrow \Abs{f(x)}<\varepsilon.
$$
Moreover, $A$ is attractive a.s. so that there exists $N_x\subset \Omega$ with $\PP\left(N_x\right) = 0$ such that, for all $p\in N_x^c$ and $ \widetilde{\varepsilon} > 0$, there exists $T_{p,\widetilde{\varepsilon}} > 0$ with $d\left(X_t^x(p),A\right) < \widetilde{\varepsilon}$ for $t>T_{p,\widetilde{\varepsilon}}$. Thus, for $\widetilde{\varepsilon} = \delta_\varepsilon$, we have that, for all $p\in N_x^c$, there exists $T_{p,\delta_\varepsilon} > 0$ such that $d(X_t^x(p),A)<\delta_\varepsilon$ for $t>T_{p,\widetilde{\varepsilon}}$. This implies that, for all $p\in N_x^c$, we have $\Abs{f\Par{X_t^x(p)}}<\varepsilon$ for $t>T_{p,\widetilde{\varepsilon}}$. Since $\varepsilon>0$ is arbitrary, this concludes the proof.  

\item \textit{Sufficiency}. The proof follows similar lines as in 2). \\
\textit{Necessity}. Let fix $f\in\mathcal{F}_{A_c}, \varepsilon\in(0,1)$ and $r>0$. Since $f$ is continuous and has zero value on $A$, there exists $\delta_{r}>0$ such that 
$$
d(x,A)<\delta_{r}\Longrightarrow \Abs{f(x)}<r/2.
$$
Moreover, $A$ is stable in probability so that there exists $\delta_{\varepsilon,\delta_{r}} > 0$ such that 
$$
d(x,A) < \delta_{\varepsilon,\delta_{r}} \Longrightarrow
\PP\left(\sup_{t\geq0} d\left(X_t^x,A\right)<\delta_r\right)\geq 1-\varepsilon.
$$
We define the sets
$$A_x^d := \left\{p \in \Omega \mid \displaystyle\sup_{t\geq0}d(X_t^x(p),A)<\delta_r\right\}$$
and 
$$
A_x^f := \left\{p\in \Omega \mid \sup_{t\geq0}\Abs{f\Par{X_t^x(p)}} <  r \right\}.
$$
Then, for all $p\in A_x^d$ and $t \geq 0$, we have $\Abs{f\Par{X_t^x(p)}} < r/2$ so that $A_x^d \subset A_x^f$ and $1-\varepsilon \leq
\PP\left(A_x^f\right),
$
which concludes the proof. 

\item \textit{Sufficiency}. The proof follows similar lines as in 2).\\
\textit{Necessity}. Let fix $f\in\mathcal{F}_{A_c}$ and $\varepsilon \in (0,1)$. For any $x\in \X$, we define the sets 
$$
B_x^d := \left\{p \in \Omega \mid \lim_{t\rightarrow\infty}d(X_t^x(p),A)=0\right\}$$
and
$$B_x^f := \left\{p \in \Omega \mid \lim_{t\rightarrow\infty}f(X_t^x(p))=0\right\}. 
$$
Then, we can show that $B_x^d \subset B_x^f$ by using the continuity of $f$ as detailed in 2). Moreover, by assumption, there exists $\delta_\varepsilon>0$ such that 
$$
d(x,A) < \delta_\varepsilon \Longrightarrow \PP\left(B_x^d\right)\geq 1-\varepsilon,
$$
which concludes the proof since $\PP\left(B_x^f\right) \geq \PP\left(B_x^d\right) \geq 1-\varepsilon$.
\end{enumerate}
\end{proof}

\section{Koopman eigenfunctions and stability} \label{sec:res_eigen}

\noindent In the previous section, we showed that the stability of a set $A\subset \X$ is characterized by the evolution of all functions $\mathcal{F}_{A_c}$ under the action of the (pathwise) Koopman operator. However, this condition might not be easy to check in practice. Instead, we can exploit the spectral properties of the operator and show that specific eigenfunctions provide substantial information on the stability of SDEs. This will allow us to restrict our analysis to these functions, rather than working on the entire subspace $\mathcal{F}_{A_c}$.

\begin{defn}
The function $\FP\in\mathcal{F}$ is an eigenfunction of the Koopman operator if and only if $K_t\phi_\lambda = e^{\lambda t}\FP$. 
\end{defn}

The Koopman generator $\A$ and semigroup $K_t$ share the same eigenfunctions $\phi_{\lambda}\in\mathcal{F}$ associated with the eigenvalues $\lambda$ and $e^{\lambda t}$, respectively, i.e. $\A\phi_{\lambda} = \lambda \phi_{\lambda}$ and $\K_t\phi_\lambda = e^{\lambda t}\phi_\lambda$. Working with $\A$ is not very convenient as $K_t$ is usually unknown. On the other hand, $\gen$ can be computed easily since $F$ and $G$ are known. Thus, since $\A f = \gen f$ for all $f\in C^2_b(\X)\subset C^2(\X)$, we will rather work with $\gen$. Moreover, we will investigate the spectral properties of $\gen$ on the less restrictive space $C^2(\X)$. In \cite{mauroy2014global}, the authors show that the zero-level sets of specific eigenfunctions are directly related to the stability of the system.
The following theorem is the stochastic counterpart of this result.  

\begin{thm}\label{thm:M0_SDE}
Assume that $\mathcal{L}$ admits $m$ real eigenvalues $\lambda_i < 0$ and the associated eigenfunctions $\phi_{\lambda_i}\in C^2(\X)$ which extend continuously to $\overline{\X}$ and are nonnegative on $\X$. If $\X$ is bounded, or otherwise if $V:=\sum_{i=1}^m\phi_{\lambda_i}$ is radially unbounded, then, the set
$$
M = \displaystyle\bigcap_{i=1}^m\left\{x\in \overline{\X}\mid \phi_{\lambda_i}(x) = 0\right\}=\left\{x\in \overline{\X}\mid V(x) = 0\right\}
$$
is almost surely forward-invariant and almost surely attractive in $\mathbb{X}$.

\end{thm}
\begin{proof}
\textit{Invariance}. Let $x\in M$. For all $i=1,\ldots,m$, we have that $\E\left[\FPi\Par{X_t^x}\right] = e^{\lambda_i t}\FPi(x) = 0$. Since $\FPi$ is positive on $\X$, $\FPi\Par{X_t^x} = 0$ a.s. This implies that for all $i$, there exists $N_{i,x}\subset\Omega$ with $\mathbb{P}(N_{i,x}) = 1$ such that
$
\FPi\left(X_t^x(p)\right)=0
$
for any $p\in N_{i,x}$. Let $N_x=\cap_{i=1}^m N_{i,x}$. Since the intersection is finite, we have $\mathbb{P}[N_x]=1$ and, for all $p\in N_x$, $V(X^x_t(p))= \sum_{i=1}^m \phi_{\lambda_i}(X_t^x(p))=0$ which implies that $X_t^x(p)\in M$.\\
\textit{Attractivity}
The proof of this result is inspired by \cite{mao1999stochastic}. For any $i=1,\cdots,m$, we have $\FPi\in C^2(\X)$ and $\{X_t^x\}_t$ exists so that for any $t\geq 0$ and $x\in \X$, Ito's lemma yields
$$
d\FPi\Par{X_t^x} = \underbrace{\lambda_i\FPi\Par{X_t^x}}_{:=\mu_{i,t}} dt + \underbrace{\left(\nabla_x^\top\FPi G\right)\Par{X_t^x}}_{:=\sigma_{i,t}}dW_t.
$$
Let us define $Y_{i,t} = \FPi\Par{X_t^x}$ and $g_i(t,y) = e^{-\lambda_i t}y$. Clearly, $g_i$ is twice continuously differentiable in $y$ and continuously differentiable in $t$. Thus, we can apply Ito's lemma so that 
$$
dg_i(t,Y_{i,t}) = e^{-\lambda_i t}\sigma_{i,t} \, dW_t,
$$
or, equivalently, 
$$
g(t,Y_{i,t}) = g(0,Y_{i,0}) + \int_0^te^{-\lambda_i s}\sigma_{i,s}dW_s.
$$
Considering that $g_i(t,y) = e^{-\lambda_i t}y$, we have
$$
e^{-\lambda_i t}\FPi\left(X_t^x\right) = \FPi(x) + \underbrace{\displaystyle\int_0^t e^{-\lambda_i s}\nabla \FPi(X_s)^\top G(X_s)\,\mathrm{d}W_s}_{:= M_{i,t}},  
$$
where $M_{i,t}$ is a local martingale satisfying $M_{i,0} = 0$. Thus, it follows from \cite[Lemma 2.2]{mao1999stochastic} that there exists $Y_{i,\infty}\geq 0$ such that 
$$
\lim_{t\rightarrow\infty}e^{-\lambda_i t}\FPi\left(X_t^x\right) = Y_{i,\infty}\quad\text{a.s.},
$$ 
so that $\FPi\left(X_t^x\right) \rightarrow 0$ almost surely. This implies that there exists $N_{i,x}\subset\Omega$ with $\mathbb{P}(N_{i,x}) = 1$ such that
$$
\lim_{t\rightarrow\infty} \FPi\left(X_t^x(p)\right)=0
$$
for any $p\in N_{i,x}$. Let $N_x=\cap_{i=1}^m N_{i,x}$. Since the intersection is finite, we have $\mathbb{P}[N_x]=1$ and, for any $p\in N_x$, it follows that
\begin{equation}\label{proof:thm2:lim_phi}
\lim_{t\to\infty} V(X^x_t(p))=\lim_{t\to\infty} \sum_{i=1}^m \phi_{\lambda_i}(X_t^x(p))=0.
\end{equation}
Let fix $p\in N_x$. We claim that $\{X_t^x(p) \mid t\geq 0\}$ must be bounded. If $\X$ is bounded, the result is straightforward since $\X$ is forward invariant. Suppose that $\X$ is unbounded and assume by contradiction that $\{X_t^x(p) \mid t\geq 0\}$ is unbounded. Then, there exists a subsequence $(t_k)_{k\in \mathbb{N}}$ such that 
\begin{equation}\label{lim_norm}
\lim_{k\rightarrow\infty}\|X_{t_k}^x(p)\|=\infty.
\end{equation}
Clearly, we have $t_k\rightarrow\infty$. Otherwise, there exist $M>0$ and a subsequence $\left(t_{k_n}\right)_{n\in \mathbb{N}}$ such that $t_{k_n}<M$ for all $n$. However, since $X_{t}^x(p)$ is continuous with respect to $t$, it is bounded on $[0,M]$. Thus, we obtain
$$\lim_{n\rightarrow\infty}\|X_{t_{k_n}}^x(p)\|\leq M,$$
which contradicts \eqref{lim_norm}. By assumption, $V$ is radially unbounded so that
$$
\lim_{t_k\to\infty}V(X_{t_k}^x(p))=\infty,
$$
which contradicts \eqref{proof:thm2:lim_phi}. Since $\{X_t^x(p) \mid t\geq 0\}$ is bounded, there is a subsequence $\{X_{t_k}^x(p) \mid k\geq 1\}$ which converges to some $y\in \overline{\X}$. Let 
$$
\omega(x,p) = \bigcap_{T\geq0}\displaystyle\overline{\left\{X_t^x(p), t\geq T\right\}}
$$
be the $\omega-$limit set of the trajectory $\{X_t^x(p)\}_t$. It is straightforward to see that $y\in\omega(x,p)$ so that $\omega(x,p)\neq\emptyset$. We now show that $\omega(x,p)\subset M_0$. Let $y \in \omega(x,p)$. Since $\FP$ is continuous, we have $\FP(y) = 0$ and $y\in M_0$. Thus, for any $p\in N_x$, we have 
$$
d(X_t^x(p),M_0)\leq d(X_t^x(p),\omega(x,p)).$$
Finally, it follows from \cite[Lemma 4.1]{KhalilNonlinearControl} that $d(X_t^x(p),\omega(x,p)) \rightarrow 0$, which implies that 
$$
\lim_{t\rightarrow\infty}d(X_t^x(p),M_0) = 0.
$$
This concludes the proof.
\end{proof}
We note that the function $V$ plays the role of a Lyapunov function. Theorem \ref{thm:M0_SDE} can also be applied to each eigenfunction $\phi_{\lambda_i}$ separately. However, each eigenfunction defined on an unbounded set $\mathbb{X}$ must be radially unbounded in this latter case, which is a condition stronger than the radial unboundedness $V$.

The existence of positive eigenfunctions required in Theorem \ref{thm:M0_SDE}  seems restrictive compared to the deterministic case, where this condition is not needed (see \cite{mauroy2014global}). In fact, positivity plays a key role in characterizing stability in both deterministic and stochastic cases. However, for deterministic systems, if $\phi_\lambda$ is an eigenfunction associated with the eigenvalue $\lambda\in \mathbb{C}$, then $|\phi_\lambda|$ is another (positive) eigenfunction associated with the eigenvalue $\Re(\lambda)$, which can be used to capture the stability properties. In contrast, for stochastic systems, the modulus of an eigenfunction is typically not an eigenfunction, so that positivity must be imposed. 

\section{Numerical examples}\label{sec:numerical}

\noindent In this section, we illustrate Theorem \ref{thm:M0_SDE} with two examples.

\begin{ex}
Consider the one-dimensional SDE
\begin{equation}\label{eq:ex1}
dX_t = -\dfrac{1}{2}\sin(2X_t) \, dt + \sqrt{2}\sin(X_t) \, dW_t.
\end{equation}
Both functions $F(x)=-\sin(2x)/2$ and $G(x)=\sqrt{2}\sin(x)$ are bounded and globally Lipschitz over $\R$. Moreover, the monotonicity condition holds since 
$$
|x F(x)| + \dfrac{1}{2}G(x)^2\leq 2x^2+x^2 \leq 3(x^2+1),
$$
so that a unique solution exists for all $t\geq0$. It can be easily verified that $\phi_\lambda(x) = \sin(x)$ is an eigenfunction of $\gen$ with eigenvalue $\lambda = -1$. Let $\X=[0,\pi]\subset\R$. Since $F(0)=G(0)=0$ and
$F(\pi)=G(\pi)=0$, a solution that reaches $0$ or $\pi$ stays on the equilibrium, so that $\X$ is forward-invariant almost surely by continuity of the trajectory. Finally, since $\phi_\lambda$ is nonnegative on $\X$ and $\X$ is
bounded, Theorem~\ref{thm:M0_SDE} applies, with
$$
M_0=\{x\in[0,\pi]\mid\sin x=0\}=\{0,\pi\}.
$$
It follows that $\lim_{t\to\infty}d\big(X_t^x,\{0,\pi\}\big)=0$ almost surely for all $x\in[0,\pi]$. By continuity of trajectories, $X_t^x$ converges almost surely either to $0$ or to $\pi$. This result is supported by numerical experiments shown in Figure \ref{fig:example2} where a realization of the process $X_t^x$ is computed for eight initial conditions $x$ sampled in $[0,\pi]$ (gray dots), using the Euler-Maruyama method. As expected, each realization converges either to $0$ or $\pi$. 

\begin{figure}[h]
\centering
\includegraphics[width=0.45\textwidth]{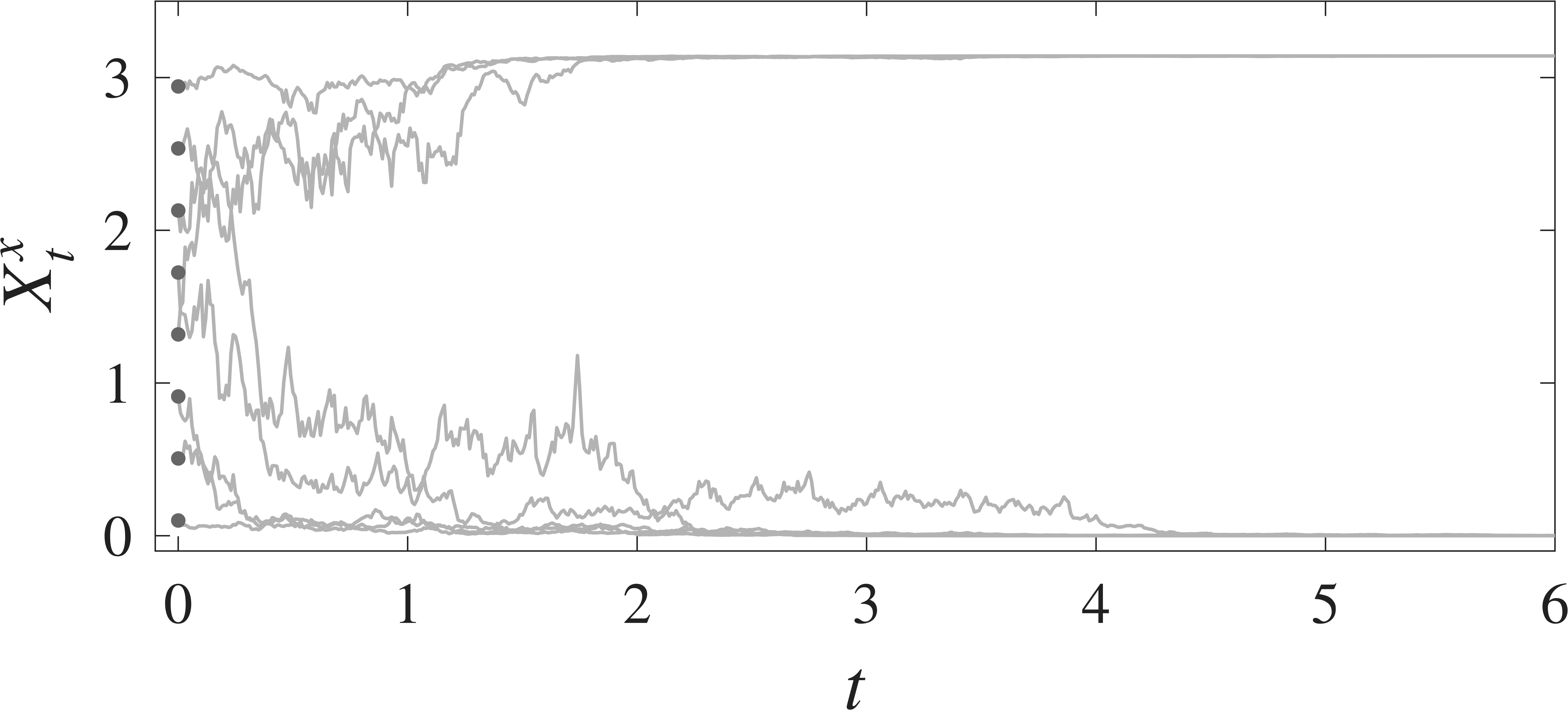}
\caption{The realizations of the solution $\{X_t^x\}$ to \eqref{eq:ex1}, computed through the Euler-Maruyama method for eight different initial conditions (gray dots), converge to $0$ or $\pi$.} 
\label{fig:example1}
\end{figure}

\end{ex}

\begin{ex}\label{ex3}

Let $a,b>0$ and consider the two-dimensional SDE \eqref{SDE} with
$$
F(x,y) = \begin{pmatrix}
-ax-\dfrac{xy^2}{2}\\[0.3cm]
-by-\dfrac{x^2y}{2}
\end{pmatrix}~~\text{and}\quad G(x,y) = \begin{pmatrix}
xy & 0\\[0.15cm]
0 & xy
\end{pmatrix}.
$$
Since both $F$ and $G$ are polynomial, the local Lipschitz condition holds. Moreover, one can easily verify that 
$$
x^\top F(x) + \dfrac{1}{2}\mathrm{Tr}(G(x)^\top G(x)) \leq 0
$$
so that the monotonicity condition holds for any $C>0$. It follows that there exists a unique solution to the SDE for all $t\geq0$. It can be easily verified that
$\phi_{\lambda_1}(x,y)=x^2$ and $\phi_{\lambda_2}(x,y)=y^2$ are the eigenfunctions of $\gen$ associated with the eigenvalues $\lambda_1=-2a$ and $\lambda_2=-2b,$ respectively.

Let $O_x=\{(x,0)\mid x\in\R\}$ and $O_y=\{(0,y)\mid y\in\R\}$. It is clear that $G$
vanishes on $O_x$, and the deterministic process $\big(xe^{-at},0\big)$ is the unique solution to the SDE
with the initial condition $(x,0)$, so that $O_x$ is forward-invariant. Similarly, $O_y$ is forward-invariant. The eigenfunction $\phi_{\lambda_1}$ is nonnegative and radially unbounded on $O_x$, and vanishes on $O_x$ only at the origin. With $\X=O_x$, Theorem \ref{thm:M0_SDE} therefore implies that the origin is almost surely attractive in $O_x$. A similar reasoning shows that the origin is almost surely attractive in $O_y$.

Theorem \ref{thm:M0_SDE} also holds on $\X = \R^2$ since $V=\phi_{\lambda_1}+\phi_{\lambda_2}$ is radially unbounded
on $\X = \R^2$. It follows that \vspace{-0.2cm}
\begin{align*}
M=O_y\cap O_x=\{(0,0)\}
\end{align*}
\vspace{-0.2cm}

\noindent is almost surely attractive in $\R^2$. This result is supported by numerical experiments shown in Figure \ref{fig:example2} where a realization of the process $X_t^x$ is computed using the Euler-Maruyama method for 10 initial conditions $x$ randomly distributed in $[-1,1]^2$. As expected, each realization converges to $0$. 

\begin{figure}[h]
\centering
\hspace{-0.4cm}\includegraphics[width=0.3\textwidth]{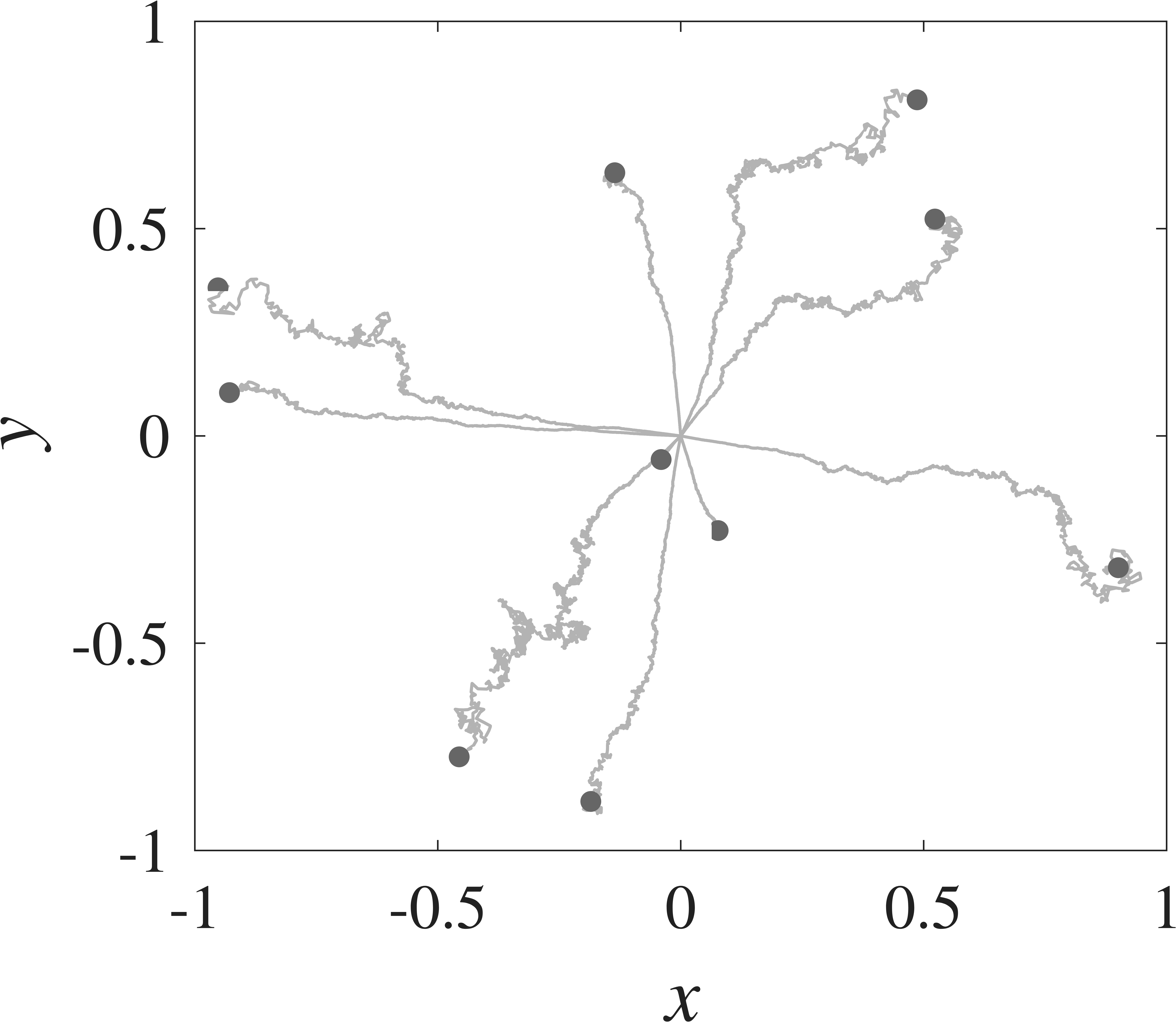}
\caption{The realizations of the solution $\{X_t^x\}$ to the SDE of Example \ref{ex3}, computed via the the Euler-Maruyama method for ten different initial conditions (gray dots), converge to the origin.}
\label{fig:example2}
\end{figure}

\end{ex}

\section{Conclusion}\label{sec:ccl}

In this paper, we derived stability results for stochastic differential equations through the Koopman operator framework. We first characterized the stability of stochastic systems via the convergence properties of the Koopman semigroup. In view of the limitations of our first approach, we subsequently investigated the spectral properties of the Koopman operator. More precisely, we established sufficient conditions under which the intersection of the zero-level sets of real non-negative Koopman eigenfunctions associated with strictly negative eigenvalues is almost surely attractive. Finally, this theoretical result was illustrated with two examples.

This work opens up several perspectives. First, both examples presented in Section \ref{sec:numerical} suggest that the stability of the set containing all the equilibria of a $n-$dimensional system could be characterized by $n$ linearly independent eigenfunctions, similarly to \cite[Proposition 2]{mauroy2014global} in the deterministic case. This property could be further investigated. Moreover, since the eigenfunctions of the Koopman operator are usually unknown, they could be computed with numerical methods relying either on a Taylor-based expansion as in \cite{mauroy2014global}, or on recent techniques such as those developed in \cite{hamzi2026computation,xu2025data}.  

\section*{AI Statement}
AI tools were used by the authors for (i) initial learning about stochastic differential equations, (ii) identifying existing work related to the stability analysis of SDEs for the Introduction, and (iii) assisting with the proof that \(M_{i,t}\) is a local martingale in Theorem \ref{thm:M0_SDE}. The remaining arguments in the proofs of Theorem \ref{thm:M0_SDE} and Proposition \ref{prop:1}, as well as the examples provided in Section \ref{sec:numerical}, are the authors' own work.

\bibliographystyle{IEEEtran}
\bibliography{ref}

\end{document}